\documentstyle[amsfonts,amssymb]{article}

\newcommand{\forces}{\Vdash}

\newcommand{\coll}{{\rm Coll}}

\newcommand{\V}{{\bf V}}

\newcommand{\lesdot}{\mathrel{\mathord{<}\!\!\raise
0.8 pt\hbox{$\scriptstyle\circ$}}}

\newcommand{\con}{{\frak c}}

\newcommand{\nwd}{{\rm nwd}}
\newcommand{\meager}{{\rm meager}}
\newcommand{\nmeager}{\neg{\rm meager}}
\newcommand{\can}{2^{\textstyle \omega}}
\newcommand{\fs}{2^{\textstyle <\!\omega}}
\newcommand{\baire}{\omega^{\textstyle \omega}}

\newcommand{\fseo}{\omega^{\textstyle <\!\omega}}

\newcommand{\Tr}{{\bf Tr}}
\newcommand{\Trw}{{\bf Trw}}

\newcommand{\conc}{{}^\frown\!}

\newcommand{\rest}{{\mathord{\restriction}}}

\newcommand{\dom}{{\rm dom}}

\newcommand{\C}{{\Bbb C}}

\renewcommand{\H}{{\cal H}}

\renewcommand{\L}{{\bf L}}

\newcommand{\p}{{\Bbb P}}

\newcommand{\q}{{\Bbb Q}}

\newcommand{\QED}{\hfill\vrule width 6pt height 6pt depth 0pt
\vspace{0.1in}}

\newcommand{\Proof}{\noindent{\sc Proof} \hspace{0.2in}}
\newtheorem{theorem}{Theorem}[section]
\newtheorem{claim}{Claim}[theorem]
\newtheorem{lemma}[theorem]{Lemma}
\newtheorem{proposition}[theorem]{Proposition}
\newtheorem{corollary}[theorem]{Corollary}
\newtheorem{conclusion}[theorem]{Conclusion}
\newtheorem{problem}[theorem]{Problem}
\newtheorem{definition}[theorem]{Definition}

\newtheorem{pwa}[theorem]{Problems we address}
\title{Ideals without ccc}
\author{
{\bf Marek Balcerzak}\thanks{\ \ \ Research partially supported by KBN grant
No PB 691/2/91}\\
Institute of Mathematics\\
L\'od\'z Technical University\\
90-924 L\'od\'z, Poland\\
\and
{\bf Andrzej Ros{\l}anowski}\thanks{\ \ \ Research partially supported by
KBN grant 1065/P3/93/04}\\
Institute of Mathematics\\
The Hebrew University of Jerusalem\\
Jerusalem, Israel\\
and\\
Mathematical Institute of Wroclaw University\\
50384 Wroclaw, Poland
\and
{\bf Saharon Shelah}\thanks{\ \ \ Research supported by  ``Basic Research
Foundation'' of The Israel Academy of Sciences and Humanities. Publication
number 512} 
\\
Institute of Mathematics\\
The Hebrew University of Jerusalem\\
Jerusalem, Israel\\
and\\
Department of Mathematics, Rutgers University\\
New Brunswick, NJ 08854, USA
}
\date\today

\begin{document}
\maketitle
\eject

\begin{abstract}
Let $I$ be an ideal of subsets of a Polish space X, containing all
singletons and possessing a Borel basis. Assuming that $I$ does not satisfy
ccc, we consider the following conditions (B), (M) and (D).  Condition (B)
states that there is a disjoint family $F\subseteq P(X)$ of size $\con$,
consisting of Borel sets which are not in $I$. Condition (M) states that
there is a Borel function $f:X\rightarrow X$ with $f^{-1}[\{ x\} ]\notin I$
for each $x\in X$. Provided that X is a group and $I$ is invariant, condition
(D) states that there exist a Borel set $B\notin I$ and a perfect set
$P\subseteq X$ for which the family $\{ B+x: x\in P\} $ is disjoint. The aim
of the paper is to study whether the reverse implications in the chain
$(D)\Rightarrow (M)\Rightarrow (B)\Rightarrow \mbox{not-ccc}$ can hold.  We
build a $\sigma $-ideal on the Cantor group witnessing $(M)\& \neg (D)$
(Section 2). A modified version of that $\sigma $-ideal contains the whole
space (Section 3). Some consistency results on deriving (M) from (B) for
"nicely" defined ideals are established (Sections 4 and 5). We show that both
ccc and (M) can fail (Theorems \ref{easynotcc} and \ref{connotcc}). Finally,
some sharp versions of (M) for invariant ideals on Polish groups are
investigated (Section 6). 
\end{abstract}

\section{Introduction}
An ideal on a space $X$ is a family $I$ of subsets of $X$ closed under
finite unions and subsets (i.e. $A,B\in I\ \Rightarrow\ A\cup B\in I$ and
$A\subseteq B, B\in I\ \Rightarrow\ A\in I$); $\sigma$-ideals are closed
under countable unions. All ideals we consider are assumed to be non
trivial, they do not contain the whole space $X$. Moreover we want them to
contain all singletons $\{x\}$ ($x\in X$). The ideal $I$ on a Polish space
$X$ is called Borel if it has a Borel basis (i.e. if for each set $A\in I$
there is a Borel subset $B$ of $X$ such that $A\subseteq B$ and $B\in I$).

The most popular Borel $\sigma$-ideals (e.g. the ideal of meager sets or the
ideal of Lebesgue null sets) satisfy the countable chain condition (ccc).
This condition says that the quotient Boolean algebra of Borel subsets of
the space modulo the $\sigma$-ideal is ccc (i.e. every family of disjoint
Borel sets which do not belong to the ideal is countable). In this paper we
are interested in ideals which do not satisfy this condition. The question
that arises here is what can be the reasons for failing ccc. The properties
(M) and (D) defined below imply that the ideal does not satisfy ccc (and
actually even more, see (B) below).
\begin{definition}
Let $I$ be an ideal on an uncountable Polish space $X$.
\begin{enumerate}
\item  We say that $I$ has {\em property (M)} if and only if\\
there is a Borel measurable function $f:X\rightarrow X$ with $f^{-1}[\{ x\}
]\notin I$ for each $x\in X$.
\item  Provided that $X$ is a Polish Abelian group and $I$ is invariant
(i.e. $A\in I$ and $x\in X$ imply $A+x\stackrel{\rm def}{=}\{ a+x:a\in A\} \in
I$),\\  
we say that $I$ has {\em property (D)} if and only if\\
there are a Borel set $B\notin I$ and a perfect (non-void) set $P\subseteq X$
such that $(B+x)\cap (B+y)=\emptyset$ for any distinct $x,y \in P$.
\end{enumerate}
\end{definition}
Properties (M) and (D) were introduced and investigated in \cite{Ba}. It was
observed that (D)$\Rightarrow$(M), if $I$ is invariant in the group $X$. Of
course, (M) implies the following condition:
\begin{description}
\item[(B)]\ \ \  there is a family $F\subseteq P(X)$ of cardinality $\con$
(the size of the continuum) of pairwise disjoint Borel subsets of $X$ that
are not in $I$.
\end{description}
In \cite {Ba} Fremlin's theorem stating the consistency of $\neg ($(B)$
\Rightarrow$(M)) is shown. However, it is unclear how his proof could be
applied  to the invariant case. The following questions arise (3 and 4 are
posed in \cite {Ba}):
\begin{pwa}
\label{pwa}
\hspace{0.5in}
\begin{enumerate}
\item Suppose $I$ is a Borel (invariant) ideal on a Polish space
(group) X, for which the ccc fails. Does $I$ satisfy (B)?
\item {\em (Remains open)}\ \ \ Is $\neg ($(B)$\Rightarrow$(M)) consistent,
for some invariant ideal ($\sigma$-ideal)?
\item {\em (Remains open)}\ \ \ Is $\neg ($(B)$\Rightarrow$(M)) provable in
ZFC, for some ideal ($\sigma$-ideal)? 
\item Does (M)$\Rightarrow$(D) hold for every invariant ideal ($\sigma$-ideal)?
\end{enumerate}
\end{pwa}
The present paper considers these questions.  We mostly restrict ourselves
to the Cantor group $\can$ with the coordinatewise addition modulo 2
(denoted further by $\oplus $, or simply, by $+$).

At first, let us show that question 1 of \ref{pwa} can have the negative
answer, if we do not require any additional properties of a $\sigma$-ideal
$I$. 

\begin{theorem}
\label{easynotcc}
\hspace{0.15in}For each cardinal $\kappa$, $\omega < \kappa <\con $, there
exists a Borel $\sigma$-ideal $I$ on $\can$ such that $\kappa$ is the
maximal cardinal for which one can find a disjoint family of size $\kappa$
of Borel sets in $\can$ that are not in $I$. Consequently, $I$ does not
satisfy both ccc and (B), and it satisfies $\kappa ^{+}$-cc.
\end{theorem}

\Proof Pick pairwise disjoint nonempty perfect sets $P_{\alpha }\subseteq
\can$, $\alpha <\kappa$. Define $I$ as follows: 
\begin{quotation}
\noindent $E\subseteq \can$ belongs to $I$ if and only if

\noindent there is a Borel set $B\subseteq\can$ such that $E\subseteq B$ and
for each $\alpha<\kappa$ the intersection $B\cap P_{\alpha}$ is meager in
$P_{\alpha}$.
\end{quotation}
Obviously, $I$ is a Borel $\sigma $-ideal on $\can$ and it does not satisfy
ccc since each set $P_{\alpha }$ is not in $I$.  Suppose that $F$ is a family
of pairwise disjoint $I$-positive Borel sets and $|F|=\kappa^{+}$. Then there
is an $\alpha <\kappa $ such that
\[|\{E\in F: E\cap P_{\alpha }\ \mbox { is non-meager in }\ P_{\alpha }\} |
=\kappa ^{+}.\]
This is impossible since the ideal of meager sets in $P_{\alpha }$ satisfies
ccc. \QED

A similar idea but in a much more special form will be used to produce the
negative answer of a modified version of problem \ref{pwa}(1), where (B) is
replaced by (M) and $I$ is a $\sigma $-ideal on $\can $ with $\Pi ^{1}_{2}$
definition (Conclusion \ref{connotcc}). We still do not know what can
happen if the invariance of $I$ is assumed.
\medskip 

\noindent{\bf Acknowledgment:} We would like to thank the referee for
very helpful comments on the presentation of the material. 

\section{The minimal $\sigma$-ideal without property (D)}
In this section we are going to answer question 4 of \ref{pwa} in negative.
It is stated in \cite{Ba} that Bukovsk\'y has reformulated the question about 
(M)$\Rightarrow$(D) by considering the $\sigma$-ideal $I_{0}$ generated by the
family
\[\begin{array}{ll}
F_{0}=\{ B\subseteq \can:&B\mbox{ is Borel and there is a perfect set $P
\subseteq \can$ such that}\\
\ &\{ B\oplus x: x\in P\}\mbox{ is a disjoint family }\}.
  \end{array}\]
Then $I_{0}$ is the minimal invariant $\sigma$-ideal without property (D).
Observe that $I_{0}$ is not trivial since $F_{0}$ is contained in the ideal of
measure zero sets (as well as in the ideal of meager sets). 

\begin{theorem}
\label{MnotD}
\hspace{0.15in}There is a continuous function $f:\can \rightarrow \can $
such that 
\[(\forall x\in \can)(f^{-1}[\{ x\} ]\notin I_{0}).\]
Consequently, the $\sigma$-ideal $I_{0}$ has property (M) and does not have
property (D). 
\end{theorem}
The rest of this section will be devoted to the proof of the above theorem. We
will break it into several steps presented in consecutive subsections, some of
these steps may be interesting {\em per se}. For simplicity, we shall write
$+$ for the addition in $\can $; also $+$ will be used for the addition of
finite sequences of zeros and ones. 

\subsection{A combinatorial lemma} 
\label{aclsubsec}
We start with defining the function $f$ which existence is postulated in
\ref{MnotD}. Its construction is very simple and based on the following
(essentially elementary) observation.

\begin{lemma}
\label{trivial}
For each $n\in \omega$ there are $N\in\omega$ and a subset $C\subseteq
2^{\textstyle N}$ such that every $n$ translates of $C$ have non-empty
intersection, and likewise for $C'=2^{\textstyle N}\setminus C$.
\end{lemma}

\Proof Let $n\geq 1$. Since $\log_2(2^x-1)<x$ for each $x>0$, we can choose 
$\varepsilon$ such that 
\[\max_{g\in [1,2^n]}\frac{\log_2(2^g-1)}{g}<\varepsilon<1.\]
Then $2^g-1<2^{g\varepsilon}$ for every $g\in [1,2^n]$. Next, take an integer
$N\geq n$ such that
\begin{description}
\item[$\oplus$]\qquad $nN+1<2^N(1-\varepsilon)$
\end{description}
(this is possible since $1-\varepsilon>0$). We claim that this $N$ is good for
our $n$. To show that there exists a suitable set $C\subseteq 2^{\textstyle
N}$ we will estimate the number of all ``bad'' sets. Fix for a moment a
sequence $\langle s_0,\ldots,s_{n-1}\rangle\subseteq 2^{\textstyle N}$. We
want to give an upper bound for the number of all subsets $D$ of
$2^{\textstyle N}$ such that 
\begin{description}
\item[$\otimes$] either\quad $\bigcap\limits_{k<n}D+s_k=\emptyset$\quad or
\quad $\bigcap\limits_{k<n}(2^{\textstyle N}\setminus D)+s_k=\emptyset$.
\end{description}
Let $G$ be the subgroup of $(2^{\textstyle N},+)$ generated by $\{s_0,\ldots,
s_{n-1}\}$, $g=|G|$. Clearly $1\leq g\leq 2^n$ and hence, by the choice of
$\varepsilon$, $2^g-1<2^{g\varepsilon}$. Suppose that a set $D\subseteq 2^{
\textstyle N}$ is such $\bigcap\limits_{k<n}D+s_k=\emptyset$. Then for each
$s\in 2^{\textstyle N}$ there is $k<n$ such that $s_k+s\notin D$. Hence, for
every $s\in 2^{\textstyle N}$, 
\[D\cap(G+s)\neq\{s_0+s,\ldots,s_{n-1}+s\}\subseteq G+s.\]
Consequently, as $\{G+s: s\in 2^{\textstyle N}\}$ is a partition of
$2^{\textstyle N}$ into $2^N/g$ sets, the number of all $D\subseteq
2^{\textstyle N}$ satisfying the condition $\otimes$ is not greater than 
\[2\cdot(2^g-1)^{2^N/g}<2\cdot 2^{g\varepsilon 2^N/g}=2^{\varepsilon 2^N+1}.\]
There are $(2^N)^n$  sequences $\langle s_0,\ldots,s_{n-1}\rangle\subseteq
2^{\textstyle N}$ and each of them eliminates less than $2^{\varepsilon
2^N+1}$ subsets of $2^{\textstyle N}$. Hence there are at most $2^{nN}\cdot
2^{\varepsilon 2^N+1}=2^{\varepsilon 2^N+nN+1}$ ``bad'' sets $D\subseteq
2^{\textstyle N}$. By $\oplus$ we have that $\varepsilon 2^N+nN+1<2^N$ so
there is a set $C\subseteq 2^{\textstyle N}$ which is ``good''. The lemma is
proved. \QED 
\medskip

\noindent Applying Lemma \ref{trivial} inductively, choose integers $n_i$ and
sets $C_i$ (for $i\in\omega$) such that
\begin{description}
\item[($\alpha$)]\hspace{3em} $0=n_{0}<n_{1}<n_{2}<...<\omega $,\quad $C_i
\subseteq 2^{\textstyle [n_i,n_{i+1})}$,
\item[($\beta$)]\hspace{3em} if $s_{0},...,s_{n_{i}}\in 2^{\textstyle [n_{i},
n_{i+1})}$ then

both $\bigcap\limits_{k\leq n_i}C_i+s_k$ and $\bigcap\limits_{k\leq n_i} (2^{
\textstyle [n_i,n_{i+1})}\setminus C_i)+s_k$ are non-empty.
\end{description}
Next define $f:\can \rightarrow \can $ by $f(x)(i)=1$ (respectively $0$) if
$x\rest [n_i,n_{i+1})$ belongs to $C_i$ (resp. does not belong to $C_i$).

In the next steps we will show that the function $f$ is as required in
\ref{MnotD}. Since, obviously, $f$ is continuous, what we have to prove is
that for every $y\in\can$ its pre-image $f^{-1}[\{ y\} ]$ is not in the ideal
$I_{0}$. As in the proof we will use the properties of the sets $C_i$ stated
in ($\alpha$), ($\beta$) above only, it should be clear from symmetry 
considerations, that it suffices to show that the set 
\[H\stackrel{\rm def}{=}\{x\in\can: (\forall i\in\omega)(f(x)(i)=1)\}\]
is not in $I_0$. (Note that the set $H$ consists of those sequences $x$ which
satisfy $x\rest [n_i,n_{i+1})\in C_i$ for all $i$ in $\omega$.) 

\subsection{A Baire topology on $H$} 
At this step, for each sequence $\bar{P}=\langle P_n: n\in\omega\rangle$ of
perfect subsets of $\can$, we introduce a topology $\tau=\tau(\bar{P})$ on
$H$. Let the sequence $\bar{P}$ be fixed in this and the next subsections.  

Let $P_{n}^{*}=P_{n}+P_{n}=\{ x+y:x,y\in P_{n}\} $ and $T_{n}^{*}=\{ x\rest
m: x\in P_{n}^{*},m\in \omega \} $ for $n\in \omega$. Note that
$P^*_n$ is a perfect set, $T_{n}^{*}\subseteq\fs $ is a perfect tree and its
body (i.e. the set of all infinite branches through the tree) is
$[T_{n}^{*}]=P_{n}^{*}$.  

In order to define the desired topology $\tau$ we need to consider tree
orderings, say $\prec$, with domain a positive integer $n$, and compatible
with the natural ordering of integers, together with an assignment of integers
$\pi(k,\ell)$ to pairs $\{k,\ell\}\in [n]^{\textstyle 2}$, where $k$ is the
immediate predecessor of $\ell$ relative to $\prec$ (and thus, in particular,
$k<\ell$). Note that the tree ordering $\prec$ can be determined from the
mapping $\pi$ alone; $\pi$ {\em is undefined} when $k$ is not the immediate
predecessor of $\ell$. Such a mapping $\pi$ will be called {\em a tree mapping
with domain $n$}, and we shall reserve the letter $\pi$, with subscripts
and/or superscripts to denote tree mappings. 

A sequence $s\in\fs$ will be called {\em acceptable} if it belongs to the tree
of $H$ and $\dom(s)$ is some $n_i$. Thus if $s$ is acceptable and
$\dom(s)=n_i$ then for each $j<i$, the restriction $s\rest [n_j,n_{j+1})$ is in
the set $C_j$. 

Let $S$ consist of all sequences $\rho=\langle\pi,s_0,s_1,\ldots,s_{n-1}
\rangle$ where $\pi$ is a tree mapping on $n$, the $s_j$'s are acceptable with
the same domain $n_i\geq n$, and $s_k+s_\ell\in T^*_{\pi(k,\ell)}$ for all
$(k,\ell)$ such that $\pi(k,\ell)$ is defined. We also set $n=n(\rho)$ and
$i=i(\rho)$. 

\begin{lemma}
\label{getext}
Suppose that $\rho =\langle\pi,s_{0},\ldots,s_{n(\rho)-1}\rangle\in S$ and $j>
i(\rho)$. Then there are $t_{0},\ldots,t_{n(\rho )-1}\in 2^{\textstyle n_{j}}$
such that $s_{0}\vartriangleleft t_{0},\ldots,s_{n(\rho )-1}\vartriangleleft
t_{n(\rho )-1}$ and $\langle\pi ,t_{0},\ldots,t_{n(\rho )-1}\rangle\in S$
(where $s\vartriangleleft t$ means that the sequence $t$ is a proper extension
of $s$). If $j$ is sufficiently large, $t_{0},\ldots,t_{n(\rho )-1}$ can be
chosen pairwise distinct.  
\end{lemma}

\Proof Let $i=i(\rho )$, $n=n(\rho )$. Let $\prec$ be the tree ordering on $n$
determined by the tree mapping $\pi$ (so $\ell$ is the immediate
$\prec$--successor of $k$ if and only if $\pi(k,\ell)$ is defined). We will
consider the case $j=i+1$ since, for greater numbers $j$, simple induction
works.  

\noindent What we have to do is to find sequences $r^k\in C_i$ (for $k<n$)
such that 
\begin{quotation}
\noindent if $k<\ell<n$ and $\ell$ is the immediate $\prec$-successor of $k$

\noindent then $(s_k\conc r^k)+(s_\ell\conc r^\ell)\in T^*_{\pi(k,\ell)}$.
\end{quotation}
For each $k,\ell<n$ such that $k$ is the immediate $\prec$--predecessor of
$\ell$ (so in particular $k<\ell$) we choose sequences $r_{k,l}\in 2^{
\textstyle [n_{i},n_{i+1})}$ such that $(s_k+s_\ell)\conc r_{k,\ell}\in
T^{*}_{\pi(k,\ell)}$ (possible as $s_k+s_\ell\in T^{*}_{\pi (k,\ell)}$ and
$T^*_{\pi(k,\ell)}$ is a perfect tree). The sequences $r_{k,\ell}$ are our
candidates for the sums $r^k+r^\ell$: 
\begin{description}
\item[$(\circledast)$] \qquad if we decide what is $r^k$ then we will put
$r^\ell=r^k+r_{k,\ell}$. 
\end{description}
As $\prec$ is a tree ordering it follows that if we keep the above rule then
the choice of the sequence $r^0$ determines all the sequences $r^1,\ldots,
r^{n-1}$. Why? Take $k<n$. Then there exists {\em the unique} sequence
$0=k_0<k_1<\ldots<k_m=k$ such that $k_{i+1}$ is the immediate
$\prec$--successor of $k_i$ (for $i<m$) and therefore, by $(\circledast)$, 
\[r^k=r^0+r_{k_0,k_1}+r_{k_1,k_2}+\ldots+r_{k_{m-1},k_m}.\]
So choosing $r^0$ we have to take care of the demand that $r_k\in C_i$ for
{\em all} of the sequences $r^k$ (for $k<n$). Thus we have to find $r^{0}\in
C_i$ such that 
\begin{description}
\item [$(\triangle )$]\quad if $0= k_{0}<\ldots <k_{m}<n$ and $k_{i+1}$ is the
immediate $\prec$--successor of $k_i$ (for $i<m$)

then $r^{0}\in C_i+r_{k_{0},k_{1}}+r_{k_{1},k_{2}}+\ldots+r_{k_{m-1},k_{m}}$.  
\end{description}
Why can we find such an $r^0$? Each positive $k<n$ appears as the largest
element in exactly one sequence $k_{0},\ldots,k_{m}$ as in $(\triangle)$, so
we get $n-1$ translations of $C_i$ in $(\triangle )$. Thus, considering one
more trivial translation (the identity function) we have $n\leq n_{i}$ of
them.  Applying condition ($\beta$) of the choice of $C_i$, $n_{i+1}$ we may
find a suitable $r^0\in 2^{\textstyle [n_i, n_{i+1})}$.

Now, as we stated before, the choice of $r^0$ and $(\circledast)$ determine
all sequences $r^0,r^1,\ldots,r^{n-1}$. Moreover 
\begin{quotation}
\noindent if $0=k_0<\ldots<k_m=k<n$ is a sequence as in $(\triangle)$ 

\noindent then $r^k=r^0+r_{k_0,k_1}+\ldots r_{k_{m-1},k_m}\in C_i$ 
\end{quotation}
(we use the fact that the addition and the subtraction in $\fs$ coincide).

Define $t_k=s_k\conc r^k$ for $k<n$. We immediately get $\langle\pi,t_0,
\ldots,t_{n-1}\rangle\in S$.  

Finally suppose that $k,\ell<n$ are such that $k$ is the immediate
$\prec$--predecessor of $\ell$ and $\ell_0<\ell$. Take $j>i$ and $r^*_{k,
\ell},r^{**}_{k,\ell}\in 2^{\textstyle [n_i,n_j)}$ such that
\begin{quotation}
\noindent $(s_k+s_\ell)\conc r^*_{k,\ell},(s_k+s_\ell)\conc r^{**}_{k,\ell}
\in T^*_{\pi(k,\ell)}\cap 2^{\textstyle n_j}$ and

\noindent $r^*_{k,\ell}\rest n_{j-1}=r^{**}_{k,\ell}\rest n_{j-1}$ but

\noindent $r^*_{k,\ell}\rest [n_{j-1},n_j)\neq r^{**}_{k,\ell}\rest [n_{j-1},
n_j)$
\end{quotation}
(possible as $T^*_{\pi(k,\ell)}$ is perfect). If we now repeat the procedure
described earlier choosing $r^*_{k,\ell}\rest [n_m,n_{m+1})$ as $r_{k,\ell}$
at the stages $m<j-1$ then, extending the sequences from $n_{j-1}$ to $n_j$ we
may use either $r^*_{k,\ell}\rest [n_{j-1},n_j)$ or $r^{**}_{k,\ell}\rest
[n_{j-1},n_j)$. Consequently we may make sure that the respective sequence
$r^\ell$ is distinct from $r^{\ell_0}$. Repeating this for all pairs $\ell_0<
\ell<n$ we may get that all the final extensions $t_\ell$ are distinct. The
Lemma is proved. \QED 
\medskip

Note that if $\rho=\langle\pi,s_0,\ldots,s_{n-1}\rangle\in S$, $n+1\leq
n_{i(\rho)}$ and $m<\omega$ then $\langle\pi',s_0,\ldots,s_{n-1},s_0\rangle
\in S$, where $\pi'$ is such that $\pi'\rest [n]^{\textstyle 2}=\pi$,
$\pi'(0,n)=m$ and $\pi'(k,\ell)$ is undefined in all remaining cases. 
(Remember that finite sequences constantly equal to $0$ are in $T^*_m$.) Thus
we may ``extend'' each element of $S$ (to an element of $S$) getting both
longer sequences $s_i$ and the number of these sequences (i.e. $n(\rho)$)
larger. It should be remarked here that $S$ is nonempty -- it is easy to give
examples $\rho$ of elements of $S$ with $n(\rho)=1$.

For $\rho \in S$ we define the basic set $U(\rho)$ as
\[\begin{array}{ll}
\{x_{0}\in H:& (\exists x_{1},\ldots,x_{n(\rho )-1}\in H)(s_{0}
\vartriangleleft x_{0}, s_{1}\vartriangleleft x_{1},\ldots,s_{n(\rho )-1}
\vartriangleleft x_{n(\rho )-1})\\
\ & \mbox{and }\ \ (\forall j>i(\rho))(\langle\pi,x_{0}\rest n_{j},\ldots,
x_{n(\rho)-1}\rest n_{j}\rangle\in S)\}\\
\end{array}\]
where $\rho =\langle\pi ,s_{0},\ldots ,s_{n(\rho )-1}\rangle$. It follows from
\ref{getext} that each $U(\rho)$ is a non-empty subset of $H$. 

\begin{proposition}
\label{gettop}
If $\rho _{1}, \rho _{2} \in S,\; x_{0}\in U(\rho_{1})\cap U(\rho _{2})$ then,
for some $\rho \in S$, we have $x_{0}\in U(\rho )\subseteq U(\rho _{1})\cap
U(\rho _{2})$.\\
Consequently the family $\{U(\rho ):\rho\in S\} $ forms a (countable) basis of
a topology on $H$. [We will denote this topology by $\tau(\bar{P})$ or just
$\tau$ if $\bar{P}$ is understood.]
\end{proposition}

\Proof For $j=0,1$, let $\rho _{j}=\langle\pi^{j},s^{j}_{0},\ldots ,s^{j}_{n^{
j}-1}\rangle$ and let $x^{j}_{1},\ldots,x^{j}_{n^{j}-1}\in H$ witness that
$x_{0}\in U(\rho_{j})$.  We shall define
$\rho $. Put $n(\rho )=n^{0}+n^{1}-1$ and $i=i(\rho )=\max\{ i(\rho
_{0}),i(\rho _{1}),n(\rho )\} +1$. Define a partial mapping $\pi$ from $[n(
\rho )]^{\textstyle 2}$ to $\omega$ as follows
\[\begin{array}{ll}
\pi\rest [n^{0}]^{\textstyle 2}=\pi^{0},    & \ \\
\pi(n^{0}-1+k,n^{0}-1+\ell)=\pi^{1}(k,\ell) & \mbox{if }0<k<\ell<n^{1},\ \{k,
\ell\}\in\dom(\pi^1),\\
\pi(0,n^{0}-1+\ell)=\pi^{1}(0,\ell)         & \mbox{if }0<\ell<n^{1},\
\{0,\ell\}\in\dom(\pi^1),\\
\pi(k,\ell)\mbox{ is undefined}             & \mbox{in the remaining cases.}
  \end{array}\]
It should be clear that $\pi$ is a tree mapping on $n(\rho)$. Finally, put
\[\rho=\langle\pi,x_{0}\rest n_{i},x^{0}_{1}\rest n_{i},\ldots,x^{0}_{n^{0}-
1}\rest n_{i},x^{1}_{1}\rest n_{i},\ldots,x^{1}_{n^{1}-1}\rest n_{i}\rangle.\]
By the choice of $i$ and $\pi$ we easily check that $\rho\in S$ and $x_{0}\in
U(\rho )$ is witnessed by $x^{0}_{1},\ldots,x^{0}_{n^{0}-1},x^{1}_{1},\ldots
x^{1}_{n^{1}-1}$ and that $U(\rho)\subseteq U(\rho_{0})\cap U(\rho_{1})$. 

To conclude the proof of the proposition note that, since each $x_{0}\in H$ is
in some $U(\rho )$ (take $n(\rho)=1$ and $s_{0}=x_0\rest n_5$), the family
$\{U(\rho ):\rho \in S\} $ is a (countable) basis of a topology on $H$. \QED

\begin{proposition}
\label{getold}
The topology $\tau$ is stronger than the product topology of $\can$ restricted
to $H$. Consequently, all (ordinary) Borel subsets of $H$ are Borel relative
to the topology $\tau$. 
\end{proposition}

\Proof Since the sets $[t]=\{x\in\can:t\vartriangleleft x\}$ for $t\in
\bigcup\limits_{i>1}2^{\textstyle n_{i}}$ form a basis of the natural topology
in $\can$, it is enough to show that $s\in 2^{\textstyle n_{i}}$, $i>1$
implies $[s]\cap H\in\tau$. But for this observe that if $[s]\cap H\neq
\emptyset$ then $n(\rho )=1$, $i(\rho )=i$ and $s_{0}=s$ generate $\rho\in S$
such that $U(\rho )=[s]\cap H$. \QED

\begin{proposition}
\label{getbaire}
$\langle H,\tau\rangle$ is a Baire space, actually each basic $\tau $-open set
$U(\rho )$ is not $\tau$-meager.
\end{proposition}

\Proof First note that each $U(\rho)$ (for $\rho\in S$) is a projection of a
compact subset of $H^{n(\rho)}$ and hence it is a (non-empty) compact set (in
the natural topology). Suppose to the contrary that for some $\rho\in S$ we
have $U(\rho )=\bigcup\limits_{k\in\omega}N_k$, where each $N_k$ is
$\tau$-nowhere dense. Then we may inductively choose $\rho_{0},\rho_{1},
\ldots\in S$ such that $\rho_0=\rho$ and
\begin{quotation}
\noindent $U(\rho_{k+1})\cap N_k=\emptyset$,\quad and $U(\rho_{k+1})\subseteq
U(\rho_k)$ \qquad for each $k\in \omega$. 
\end{quotation}
It is possible as if $\rho^*\in S$, $N$ is $\tau$-nowhere dense then there is
a non-empty $\tau$-open set $U\subseteq U(\rho^*)\setminus N$ (remember that
$U(\rho^*)\neq\emptyset$). But the sets $U(\rho')$ for $\rho'\in S$ constitute
the basis of the topology $\tau$, so we find $\rho'\in S$ with $U(\rho')
\subseteq U\subseteq U(\rho^*)\setminus N$. 

\noindent Now, by the compactness of the sets $U(\rho_k)$ we get
\[\emptyset\neq\bigcap_{k\in\omega}U(\rho_k)\subseteq U(\rho)\setminus
\bigcup_{k\in\omega} N_k\]
-- a contradiction finishing the proof of the proposition. \QED

\subsection{$\tau$--non-meager subsets of $H$}

\begin{proposition}
\label{getmeager}
If $B\subseteq H$ is a $\tau$--non-meager set with the Baire property (with
respect to $\tau$) then for every $m\in\omega$ there are distinct $x,y\in P_m$
such that the intersection $(B+x)\cap (B+y)$ is non-empty.
\end{proposition}

\Proof Let $m\in\omega$ and $B\subseteq H$ be a $\tau$--non-meager set with
the $\tau$--Baire property. Then for some $\rho\in S$ the set $U(\rho)
\setminus B$ is $\tau$-meager. Let $U(\rho)\setminus B=\bigcup\limits_{k\in
\omega}N_k$ where each $N_k$ is $\tau$-nowhere dense. Let $\rho=\langle\pi,
s_{0},\ldots s_{n(\rho )-1}\rangle$. Put 
\[\rho '=\langle\pi ',s_{0},\ldots,s_{n(\rho )-1},s_{0},\ldots
,s_{n(\rho )-1}\rangle\]
where $\pi'$ is a tree mapping on $2n(\rho)$ given by
\[\begin{array}{ll}
\pi'\rest [n(\rho)]^{\textstyle 2}=\pi,     &\ \\
\pi'(n(\rho)+k,n(\rho)+\ell)=\pi(k,\ell)    &\mbox{ if }k<\ell<n(\rho),\ \{k,
\ell\}\in\dom(\pi),\\
\pi'(0,n(\rho))=m,                          &\ \\
\pi'(k,\ell)\mbox{ is undefined}            &\mbox{in the remaining cases.}\\
  \end{array}\]
Easily, $\pi'$ is indeed a tree mapping and $\rho'\in S$ (remember that
$s_0+s_0=\bar{0}\in T^*_m$).

For a partial function $\pi_0$ from $[n]^{\textstyle 2}$ to $\omega$,
$n>n(\rho)$ let $(\pi_0)^*$ be defined by 
\[\begin{array}{l}
\dom((\pi_0)^*)=\{\{\sigma(k),\sigma(\ell)\}:\{k,\ell\}\in\dom(\pi_0)\}\quad
\mbox{ and}\\
(\pi_0)^*(k,\ell)=\pi(\sigma(k),\sigma(\ell)),
  \end{array}\]
where $\sigma:n\longrightarrow n$ is a permutation of $n$ given by
\[\sigma(0)=n(\rho),\quad \sigma(i+1)=i\mbox{ for }0\leq i<n(\rho)\quad\mbox{
and }\sigma(i)=i\mbox{ for }n(\rho)<i<n.\]
Note that if $\pi_0$ is a tree mapping on $n>n(\rho)$ such that $\{0,
n(\rho)\}\in\dom(\pi_0)$ then $(\pi_0)^*$ is a tree mapping too. Hence, if
$\rho_0=\langle\pi_0,t_0,\ldots,t_{n-1}\rangle\in S$, $n(\rho)<n$ and $\pi_0(
0,n(\rho))$ is defined then 
\[(\rho_0)^*\stackrel{\rm def}{=}\langle (\pi_0)^*,t_{n(\rho)},t_0,\ldots,
t_{n(\rho)-1},t_{n(\rho)+1},\ldots,t_{n-1}\rangle \mbox{ is in }S\mbox{ too.}\]

\begin{claim}
\label{cl1}
Suppose that $\rho_0=\langle\pi_0,t_0,\ldots,t_{n-1}\rangle\in S$ is such that
$n>n(\rho)$ and $\pi_0(0,n(\rho))$ is defined. Let $N\subseteq H$ be a
$\tau$--nowhere dense set. Then there is $\rho^+=\langle\pi^+,s^+_0,\ldots,
s^+_{k-1}\rangle\in S$ such that
\begin{enumerate}
\item $U(\rho^+)\subseteq U(\rho_0)\setminus N$,\quad $U((\rho^+)^*)\subseteq
U((\rho_0)^*)\setminus N$,\quad and 
\item $i(\rho^+)>i(\rho_0)$, $k=n(\rho^+)>n(\rho_0)=n$, $\pi^+\rest
[n(\rho_0)]^{\textstyle 2}=\pi_0$,\quad  and
\item if $\ell<n$ then $t_\ell\vartriangleleft s^+_\ell$.
\end{enumerate}
[The operation $(\cdot)^*$ is as defined before.]
\end{claim}

\noindent {\em Proof of the claim:}\ \ \  Since $N$ is $\tau$--nowhere dense
we find $\rho_1$ in $S$ such that $U(\rho_1)\subseteq U(\rho_0)\setminus N$. 
Applying the procedure from the proof of Proposition \ref{gettop} (with an
arbitrary $x_0\in U(\rho_1)$) we may find $\rho_2=\langle\pi_2, s^2_0,\ldots,
s^2_{n(\rho_2)-1}\rangle\in S$ such that $U(\rho_2)\subseteq U(\rho_1)$ and
$i(\rho_2)>i(\rho_0)$, $n(\rho_2)>n(\rho_0)$, $\pi_2\rest [n(\rho_0)]^{
\textstyle 2}=\pi_0$, $t_\ell\vartriangleleft s^2_\ell$ for $\ell<n$. Now we
look at $(\rho_2)^*\in S$ and we choose $\rho_3\in S$ such that   
\[U(\rho_3)\subseteq U((\rho_2)^*)\setminus N\subseteq U((\rho_0)^*)\setminus
N.\] 
Next, similarly as $\rho_2$, we get $\rho_4=\langle\pi_4,s^4_0,\ldots,
s^4_{n(\rho_4)-1}\rangle\in S$ with the corresponding properties with respect
to $(\rho_2)^*$ and $\rho_3$. So, in particular, $\pi_4\rest [n(\rho_2)]^{
\textstyle 2}=(\pi_2)^*$, $s^2_\ell\vartriangleleft s^4_{\ell+1}$ for $\ell<
n(\rho)$, $s^2_{n(\rho)}\vartriangleleft s^4_0$ and $s^2_\ell\vartriangleleft
s^4_\ell$ for $n(\rho)<\ell<n(\rho_2)$. Finally we apply the inverse operation
to $(\cdot)^*$ and we get $\rho^+\in S$ such that $(\rho^+)^*=\rho_4$. (This is
possible, as the only $j\in (0,n(\rho))$ for which the value of $\pi_4(0,j)$
is defined, is $j=1$.) It should be clear that the $\rho^+$ is as required in
the claim.
\medskip

Now, by induction on $k<\omega$, we choose $n_k$, $s^k_i$ (for $i<n_k$),
$\pi_k$ and $\rho_k$ such that 
\begin{description}
\item[(i)]\ \ \ $n_0=2n(\rho)$, $s^0_i=s^0_{n(\rho)+i}=s_i$ (for
$i<n(\rho)$), $\pi_0=\pi'$, 

\noindent [So $\rho'=\langle\pi_0,s^0_0,\ldots,s^0_{n_0-1}\rangle=\rho_0$.]

\item[(ii)]\ \ $\rho_k=\langle\pi_k,s^0_k,\ldots,s^0_{n_k-1}\rangle\in S$,
$U(\rho_k)\subseteq U(\rho)$, $U((\rho_k)^*)\subseteq U(\rho)$,

\noindent [Here, $(\rho_k)^*$ is the element of $S$ obtained from $\rho_k$ by
moving $s^k_{n(\rho)}$ to the first place, see the definition of $(\pi_0)^*$,
$(\rho_0)^*$ above.]

\item[(iii)]\ $n_k< n_{k+1}$, $i(\rho_k)<i(\rho_{k+1})$,
$s^k_i\vartriangleleft s^{k+1}_i$ for $i<n_k$, $\pi_k=\pi_{k+1}\rest [n_k]^{
\textstyle 2}$,

\item[(iv)]\ \ $U(\rho_{k+1})\cap N_k=U((\rho_{k+1})^*)\cap N_k=\emptyset$
and $s_i^{k+1}$ (for $i<n_{k+1}$) are pairwise distinct. 
\end{description}
The first step of the construction is fully described in the demand {\bf (i)} 
above (note that then {\bf (ii)} is satisfied as $U((\rho_0)^*)\subseteq
U(\rho)$ and {\bf (iii)}, {\bf (iv)} are not relevant). 

\noindent Suppose that we have defined $\rho_k$ etc. Apply \ref{cl1} to
$\rho_k$, $N_k$ standing for $\rho_0$, $N$ there to get $\rho_{k+1}$
(corresponding to $\rho^+$ there). It should be clear that the requirements
{\bf (ii)}--{\bf (iv)} are satisfied except perhaps the last demand of {\bf
(iv)} -- the sequences $s^{k+1}_i$ do not have to be pairwise distinct. But
this is not a problem as by Lemma \ref{getext} we may take care of this
extending them further. 
\medskip

\noindent Let $x_i=\bigcup\limits_{k\in\omega} s^k_i$ for $i<\omega$. Then
$x_i\in H$ (by {\bf (iii)}+{\bf (ii)}), $x_0\in\bigcap\limits_{k\in\omega}
U(\rho_k)$, $x_{n(\rho)}\in\bigcap\limits_{k\in\omega} U((\rho_k)^*)$ (by the
definition of $(\rho_k)^*$). Hence, by {\bf (ii)}+{\bf (iv)}, 
\[x_0,x_{n(\rho)}\in U(\rho)\setminus\bigcup_{k\in\omega} N_k\subseteq B\]
and by the last part of {\bf (iv)} they are distinct. Since $\pi_k(0,n(\rho))
=m$ for each $k<\omega$ (by {\bf (iii)}) we may apply the definition of $S$ to
conclude that $x_0+x_{n(\rho)}\in [T^{*}_{m}]$. So, $x_0+x_{n(\rho)}=x+y$ for
some $x,y\in P_m$. As $x_0\neq x_{n(\rho)}$, also $x\neq y$. We finish the
proof of the proposition noting that $x_0+x=x_{n(\rho)}+y$ and $x_{0}+x\in
B+x$, $x_{n(\rho )}+y\in B+y$, as required. \QED 

\subsection{Conclusion of the proof of Theorem \ref{MnotD}}
As we stated at the end of the subsection \ref{aclsubsec}, to conclude Theorem
\ref{MnotD} it is enough to show that the set $H$ defined there is not in the
ideal $I_0$. Suppose to the contrary that $H$ can be covered by the union
$\bigcup\limits_{n\in\omega} B_n$ of Borel (in the standard topology) subsets
$B_n$ of $\can$, each $B_n$ from the family $F_0$. Then, for every
$n\in\omega$ we may pick up a perfect set $P_n\subseteq\can$ witnessing
``$B_n\in F_0$''. Let $\bar{P}=\langle P_n: n\in\omega\rangle$ and let $\tau=
\tau(\bar{P})$ be the topology on $H$ determined by $\bar{P}$ in Proposition
\ref{gettop}. We know that the space $H$ is not meager in this topology (by
Proposition \ref{getbaire}) and therefore one of the sets $B_n\cap H$ is
$\tau$--non-meager, say $B_{n_0}\cap H$. But $B_{n_0}\cap H$ is Borel in the
standard topology of $H$, so it has the $\tau$--Baire property. Applying
Proposition \ref{getmeager} to it we conclude that there are distinct $x,y\in
P_{n_0}$ such that the intersection $(B_{n_0}+x)\cap (B_{n_0}+y)$ is
non-empty, a contradiction to the choice of $P_{n_0}$. \QED 

\section{Non-Borel case}
Now, let us omit the assumption that the sets of the family $F_{0}$ defined in
the previous section are Borel. We thus get
\[\begin{array}{ll}
F_{0}^{*}=\{A\subseteq\can: &\mbox{there is a perfect set }P\subseteq\can
\mbox{ such that}\\
\ &\{A\oplus x: x\in P\}\mbox{ is a disjoint family }\}.
  \end{array}\]
Let $I_{0}^{*}$ be the $\sigma$-ideal generated by $F_{0}^{*}$. It turns out
that $I_{0}^{*}$ is not a proper ideal.

\begin{theorem}
$\can \in I_{0}^{*}$.
\end{theorem}

\Proof Here $+$ and $-$ will stand for the respective operations on ordinals,
while the addition and subtraction in $\can$ are denoted by $\oplus$,
$\ominus$, respectively. Let us define an increasing sequence $\langle
\gamma_\beta:\beta<\omega+\omega\rangle$ of ordinals as follows:
\begin{quotation}
$\gamma_{0}=0$, $\gamma_{n+1}=\gamma_{n}+\con$\qquad for $n\in\omega$,

$\gamma_{\omega}=\sup\limits_{n\in\omega}\gamma_{n}$, $\gamma_{\omega+n+1}=
\gamma_{\omega+n}+\con$\qquad for $n\in\omega$.
\end{quotation}
Let $P\subseteq\can$ be a perfect set independent in the Cantor group; cf.
\cite{My}. Pick pairwise disjoint perfect sets $P_n$ (for $n\in\omega$) such
that $\bigcup\limits_{n\in\omega}P_n=P$ and fix enumerations $P_{n}=
\{x_{\alpha}:\gamma_{n}\leq\alpha<\gamma_{n+1}\}$ (for $n\in\omega$; so
$x_{\alpha}$'s are distinct). Extend $P$ to $H$, a maximal independent set
called here a Hamel basis. We may assume that $|H\setminus P|=\con$. Let
$H\setminus P=\{x_{\alpha }:\gamma_{\omega}\leq\alpha<\gamma_{\omega+\omega}
\}$, where $x_{\alpha }$'s are distinct.  Consequently, $H=\{x_{\alpha}:
\alpha<\gamma_{\omega+\omega}\}$.

Now, let us define $y_{\alpha }$ for $\alpha<\gamma_{\omega+\omega}$, as
follows:

\begin{itemize}
\item if $\gamma_{\omega}\leq\alpha<\gamma_{\omega+\omega}$,\quad put
$y_{\alpha}=x_{\alpha}$, 
\item if $\gamma_{n}\leq\alpha<\gamma_{n+1}$, $n\in \omega$,\quad put

$y_\alpha=x_\alpha\oplus x_{\gamma_{\omega+n}+(\alpha-\gamma_n)\omega+1}
\oplus\ldots\oplus x_{\gamma_{\omega+n}+(\alpha-\gamma_n)\omega+n}$.
\end{itemize}
Then $\{y_\alpha:\alpha<\gamma_{\omega+\omega}\}$ forms a Hamel basis. Each
$y\in\can$ has the unique representation of the form $\sum\limits_{\alpha\in
u_y} y_\alpha$ where $u_y\subseteq\gamma_{\omega+\omega }$ is finite. For
$m\in\omega$ let
\[A_m=\{y\in\can: |u_y|=m\}.\]
Obviously, $\can=\bigcup\limits_{m\in\omega} A_m$ and
\begin{description}
\item[(1)]\quad $A_m\ominus A_m\subseteq\bigcup\limits_{i\leq 2m} A_{i}$.
\end{description}
The proof will be complete, if we show that
\begin{description}
\item[(2)]\quad $\{A_m\oplus x: x\in P_{m+1}\} $ is a disjoint family.
\end{description}
To get this, observe first that
\begin{description}
\item[(3)]\quad if $x,x'\in P_n$, $x\neq x'$ then $x\ominus x'\in A_{2n+2}$.
\end{description}
Indeed, let $x=x_\alpha$ and $x'=x_\beta$ for some $\alpha,\beta$ such that
$\gamma_n\leq\alpha<\gamma_{n+1}$, $\gamma_n\leq\beta<\gamma_{n+1}$ and
$\alpha\neq\beta$. Since 
\[\begin{array}{lcl}
x_{\alpha}&=&y_{\alpha}\ominus x_{\gamma_{\omega+n}+(\alpha-\gamma_{n})
\omega+1}\ominus\ldots\ominus x_{\gamma_{\omega+n}+(\alpha-\gamma_{n})\omega
+n}\\ 
x_{\beta}&=&y_{\beta}\ominus x_{\gamma_{\omega+n}+(\beta-\gamma_{n})\omega+1}
\ominus\ldots\ominus x_{\gamma_{\omega+n}+(\beta-\gamma_{n})\omega+n},\\
\end{array}\]
we conclude that $x\ominus x'\in A_{2n+2}$ (remember that $\beta-\gamma_n\neq
\alpha-\gamma_n$) and (3) is proved. To show (2), take distinct $x,x'\in
P_{m+1}$ and suppose that the intersection $(A_m\oplus x)\cap (A_m\oplus x')$
is non-empty. Then $a\oplus x=b\oplus x'$ for some $a,b\in A_m$. In other
words, $a\ominus b=x'\ominus x$. By (1), we have $x'\ominus x\in
\bigcup\limits_{i\leq 2m}A_{i}$ and, by (3), we get $x'\ominus x\in A_{2m+4}$,
a contradiction.\QED
\medskip

\noindent{\sc Remark:} For the Cantor group, the operations $\oplus$ and
$\ominus$ are identical.  However, we have distinguished them since the same
proof works for any Abelian Polish group admitting a perfect independent
set. Note that a number of groups different from $\can$ are good: by
\cite{Md}, each connected Abelian Polish group which has an element of
infinite order admits a perfect independent set.

\section{Getting (M) from ``not ccc''}
In this section we try to conclude (consistently) the property (M) from the
property (B) for nicely defined Borel ideals. The results here are
complementary, in a sense, to Fremlin's theorem mentioned in
Introduction. (But note that here we deal with ideals with simple definitions,
while the ideal constructed by Fremlin is very complicated.)

\begin{theorem}
\label{omega17}
Let $n\geq 2$. The following statement is consistent with ZFC
$+\con=\omega_n$:  
\begin{quotation}
\noindent {\em if} $B\subseteq\can\times\can$ is a $\Sigma^1_3$ set such that,
for some set $A\subseteq\can$ of size $\omega_2$, the
sections $B_x$, $x\in A$ are nonempty pairwise disjoint (where
$B_x=\{y\in\can: \langle x,y\rangle \in B\}$) 

\noindent{\em then} for some perfect set $P\subseteq\can$, the sections $B_x$,
$x\in P$ are nonempty pairwise disjoint.
\end{quotation}
\end{theorem}

\Proof Start with the universe $\V$ satisfying CH.\\
Let $\C_{\omega_n}$ be the (finite support) product of $\omega_n$ copies of
the Cohen forcing notion $\C$ and let $\langle c_\alpha:\alpha<\omega_n
\rangle$ be the sequence of Cohen reals, $\C_{\omega_n}$-generic over $\V$. 
Work in $\V[c_\alpha:\alpha<\omega_n]$.  

Clearly $\con=\omega_n$. Suppose that $B\subseteq\can\times\can$ is a
$\Sigma^1_3$ set and $\langle x_\zeta:\zeta<\omega_2\rangle$ is a sequence 
of reals such that $\zeta_1<\zeta_2<\omega_2$ implies that the sections
$B_{x_{\zeta_1}}$, $B_{x_{\zeta_2}}$ are nonempty disjoint. Let $U\in\V\cap
[\omega_n]^{\textstyle \omega}$ be such that the parameters of the
$\Sigma^1_3$ formula $\Phi(x,y)$ defining $B$ are in $\V[c_\alpha:\alpha\in
U]$. Next choose a sequence $\langle U_\zeta:\zeta<\omega_2\rangle\in\V$ of
countable subsets of $\omega_n\setminus U$ such that each real $x_\zeta$
is in $\V[c_\alpha: \alpha\in U\cup U_\zeta]$. Moreover, we demand that
\[\V[c_\alpha:\alpha\in U\cup U_\zeta]\models (\exists z)\Phi(x_\zeta,z)\]
(remember that $\Phi$ is $\Sigma^1_3$; use Shoenfield's absoluteness). By CH
in $\V$ and the $\Delta$-lemma we find $A\in [\omega_2]^{\textstyle\omega_2}$
and a countable set $U^*\supseteq U$ such that $\{U_\zeta\setminus U^*:\zeta
\in A\}$ is a disjoint family. We may assume that all sets $U_\zeta\setminus
U^*$ are infinite (by adding more members). Let $\V'=\V[c_\alpha:\alpha\in
U^*]$. 

Each sequence $\langle c_\alpha:\alpha\in U_\zeta\setminus U^*\rangle$ is
essentially one Cohen real; denote this real by $d_\zeta$. Note that if
$\zeta_0,\zeta_1,\zeta_2\in A$ are distinct then $\langle d_{\zeta_0},
d_{\zeta_1},d_{\zeta_2}\rangle$ is $\C\times\C\times\C$-generic over $\V'$.
A real $x\in\can$ in the one Cohen real extension is the value of a Borel
function $f:\can\rightarrow\can$ from the ground model at this Cohen real.
Consequently, we find a sequence $\langle f_\zeta:\zeta\in
A\rangle\in\V'$ of Borel functions from $\can$ into $\can$ such that
\[\V[c_\alpha:\alpha<\omega_n]\models x_\zeta=f_\zeta(d_\zeta).\]
By CH in $\V'$ we find $A^*\in [A]^{\textstyle \omega_2}\cap\V'$ and a Borel
function $f:\can\rightarrow\can$, $f\in\V'$ such that $f_\zeta=f$ for
$\zeta\in A^*$. By Shoenfield's absoluteness we have that for distinct
$\zeta_0,\zeta_1,\zeta_2\in A^*$:
\[\begin{array}{l}
\V'[d_{\zeta_0},d_{\zeta_1},d_{\zeta_2}]\models\\
\mbox{`` }\neg(\exists z)[\Phi(f(d_{\zeta_0}),z)\ \&\
\Phi(f(d_{\zeta_1}),z)]\ \&\ (\exists z)\Phi(f(d_{\zeta_0}),z)\ \&\
(\exists z)\Phi(f(d_{\zeta_1}),z)\mbox{ ''}
\end{array}\]
(remember the choice of $U$ and $U_\zeta$: the witness for $(\exists
z)(\Phi(x_\zeta,z))$ is in $\V[c_\alpha:\alpha\in U\cup U_\zeta]$ already). 

As the $d_\zeta$'s are Cohen reals over $\V'$ and their supports are disjoint,
by density argument we get that in $\V'$: 
\begin{description}
\item[$(*)$]\ \ \ $\forces_{\C\times\C\times\C}$

`` $\neg(\exists z)[\Phi(f(\dot{c}_0),z)\ \&\ \Phi(f(\dot{c}_1),z)]\ \&\
(\exists z)\Phi(f(\dot{c}_0),z)\ \&\ (\exists z)\Phi(f(\dot{c}_1),z)$ '',
\end{description}
where $\langle \dot{c}_0,\dot{c}_1,\dot{c}_2\rangle$ is the canonical
$\C\times\C\times\C$-name for the generic triple of Cohen reals. In
$\V[c_\alpha:\alpha<\omega_n]$ take a perfect set $P\subseteq\can$ such
that $P\times P\setminus \Delta\subseteq O$ for every open dense subset $O$
of $\can\times\can$ coded in $\V'$ ($\Delta$ stands for the diagonal
$\{\langle x,x\rangle :x\in\can\}$). Then $\langle x,y\rangle $ is
$\C\times\C$-generic over $\V'$ for each distinct $x,y\in P$. (Such a
perfect set is added by one Cohen real; the property that it is a perfect
set of mutually Cohen reals over $\V'$ is preserved in passing to an
extension.) We claim that for distinct $x,y\in P\cap\V[c_\alpha:\alpha<
\omega_n]$:
\[\begin{array}{ll}
\ &\V[c_\alpha:\alpha<\omega_n]\models\\
\ &\mbox{`` }(\exists t)\Phi(f(x),t)\ \&\ (\exists t)\Phi(f(y),t)\ \&\
\neg(\exists t)[\Phi(f(x),t)\ \&\ \Phi(f(y),t)]\mbox{ ''.}\\
\end{array}\]
Why? As $\langle x,y\rangle$ is $\C\times\C$-generic, by upward absoluteness
for $\Sigma^1_3$ formulas and $(*)$ we get
\[\V[c_\alpha:\alpha<\omega_n]\models (\exists t)\Phi(f(x),t)\ \&\
(\exists t)\Phi(f(y),t).\]
Assume that
\[\V[c_\alpha:\alpha<\omega_n]\models(\exists t)(\Phi(f(x),t)\ \&\
\Phi(f(y),t)).\]
The formula here is $\Sigma^1_3$, so we have a real $z\in\can\cap
\V[c_\alpha:\alpha<\omega_1]$ such that (by $\Pi^1_2$ absoluteness)
\begin{description}
\item[$(**)$] $\V'[x,y,z]\models(\exists t)[\Phi(f(x),t)\ \&\ \Phi(f(y),t)]$.
\end{description}
This real is added by one Cohen real over $\V'[x,y]$. So we may choose $z$ to
be a Cohen real over $\V'[x,y]$ ($\Sigma^1_3$-upward absoluteness again). Then
$\langle x,y,z\rangle$ is $\C\times\C\times\C$-generic over $\V'$ and $(**)$
contradicts $(*)$. 

Now, working in $\V[c_\alpha:\alpha<\omega_n]$, we see that for
distinct $x,y\in P$ the sections $B_{f(x)}$, $B_{f(y)}$ are disjoint
and nonempty. Hence the function $f$ is one-to-one on the perfect $P$
and we can easily get a required perfect $P'$. \QED

\begin{definition}
We say that a Borel ideal $I$ on $\can$ has a $\Pi^1_n$ definition if there is
a $\Pi^1_n$-formula $\Phi(x)$ such that
\smallskip

\centerline{$\Phi(a)\equiv$ \ ``$a$ is a Borel code and the set $\# a$ coded
by $a$ belongs to $I$''.}
\smallskip

\noindent We say that $I$ has a projective definition if it has a $\Pi^1_n$
definition for some $n$.
\end{definition}

\begin{corollary}
\label{corom}
Let $n\geq 2$. It is consistent with ZFC $+\con=\omega_n$ that:
\begin{quotation}
\noindent for each Borel ideal $I$ on $\can$ with a $\Pi^1_3$ definition, if
there exists a sequence $\langle B_\alpha: \alpha<\omega_2\rangle$ of
disjoint Borel sets not belonging to $I$ then $I$ satisfies (M).
\end{quotation}
In particular the statement 
\begin{center}
``(B) $\ \Rightarrow\ $ (M) for Borel ideals which are $\Pi^1_3$'' 
\end{center}
is consistent.
\end{corollary}

\Proof Work in the model of \ref{omega17} (so after adding $\omega_n$ Cohen
reals to a model of CH). Suppose that $I$ is a Borel ideal with $\Pi^1_3$
definition and let $\Psi(x)$ be the $\Pi^1_3$ formula witnessing it. Let

\centerline{$\Phi(x,y)\equiv$``$x$ is a Borel code\ \&\ $\neg\Psi(x)\
\&\ y\in\# x$''.}

\noindent This is a $\Sigma^1_3$ formula defining a $\Sigma^1_3$ subset $B$ of
$\can\times\can$. If there are $\omega_2$ pairwise disjoint $I$-positive Borel
sets then they determine $\omega_2$ pairwise disjoint nonempty sections of the
set $B$. By \ref{omega17} we can find a perfect set
$P\subseteq\can$ such that $\{B_x:x\in P\}$ is a family of disjoint nonempty
sets. Define a function $f:\can\rightarrow P$ by:
\begin{quotation}
\noindent for $y\in\can$, if there is $x\in P$ such that $(x,y)\in B$
then $f(y)$ is this (unique) $x$, otherwise $f(y)=x_0$ where $x_0$
is a fixed element of $P$.
\end{quotation}
Note that the set $\{(x,y)\in B: x\in P\}$ is Borel and has the property that
its projection onto $y$'s axes is one-to-one. Consequently the set
\[\{y\in\can: (\exists x\in P)((x,y)\in B)\}\]
is Borel and hence the function $f$ is Borel. Thus $f$ witnesses (M) for the
ideal $I$. (Note that no harm is done that the function is onto $P$ instead of
$\can$.) \QED

\noindent{\sc Remark:}\ \ \ 1) If we start with a model for CH and add
simultaneously $\omega_n$ random reals over $\V$ (by the measure algebra on
$2^{\omega_n}$) then in the resulting model we will have a corresponding
property for $\Sigma^1_2$ subsets of the plane $\can\times\can$ and the ideals
with $\Pi^1_2$ definitions. The proof is essentially the same as in
\ref{omega17}. The only difference is that the perfect set $P$ is a perfect
set of ``sufficiently random'' reals. For a countable elementary submodel
$N\in\V'$ of $\H(\chi)^{\V'}$ we choose by the Mycielski theorem
(cf. \cite{My2}) a perfect set $P\in\V[r_\alpha:\alpha<\omega_2]$ such that
each two distinct members of $P$ are mutually random over $N$. 

\noindent 2) Note that the demand that the ideal $I$ in \ref{corom} has to
admit $\omega_2$ disjoint Borel $I$--positive sets is not an accident. By
\ref{connotcc}, the failure of ccc is not enough.

\noindent 3) In the presence of large cardinals we may get more, see \ref{wcom}
below.

\begin{theorem}
\label{wcom}
Assume that $\kappa$ is a weakly compact cardinal and $\theta\geq\kappa^+$
is a cardinal such that $\theta^{<\kappa}=\theta$. Then there exists a
$\kappa$--cc forcing notion $\p$ which forces
\begin{quotation}
\noindent `` $\con=\theta$ and 

\noindent if $I$ is a Borel ideal with a projective definition such that there
is a sequence $\langle B_\alpha: \alpha<\omega_2\rangle$ of disjoint Borel
sets not belonging to $I$ then $I$ satisfies (M) ''
\end{quotation}
\end{theorem}

\Proof The forcing notion $\p$ is the limit $\p_\theta$ of the finite support
iteration iteration $\langle\p_0,\dot{\q}_\alpha: \alpha<\theta\rangle$ such
that $\p_0=\coll(\omega,{<}\kappa)$ and each $\dot{\q}_\alpha$ is the Cohen
forcing notion. Now repeat the arguments of \ref{omega17} with $\V^{\p_0}$ as
the ground model, to conclude that in $\V^{\p_\theta}$
\begin{quotation}
\noindent `` if $B\subseteq\can\times\can$ is a projective set which has
$\omega_2$ disjoint non-empty sections then $B$ has a perfect set of disjoint
non-empty sections. '' 
\end{quotation}
The point is that projective formulas are absolute between intermediate models
containing $\V^{\p_0}$  (see e.g. the explanations to \cite[6.5]{JuRo}). We
finish as in \ref{corom}. \QED 

\begin{theorem}
Assume that $I$ is a Borel ideal with $\Pi^1_2$ definition. Suppose
that there exists a Borel function $f:\can\rightarrow\can$ such that
$f^{-1}[\{x_\alpha\}]\notin I$ for some distinct points $x_\alpha
,\alpha<\omega_2$. Then there is a perfect set $P\subseteq\can$ such that 
\[(\forall x\in P)(f^{-1}[\{x\}]\notin I).\]
Consequently, the ideal $I$ has property (M).
\end{theorem}

\Proof Let $\Phi(x)$ be the $\Pi^1_2$ definition of the ideal $I$.
Consider the relation $E$ on $\can$:
\[x E y\ \iff\ (x=y\ \mbox { or }\ (\Phi(f^{-1}[\{x\}])\
\&\ \Phi(f^{-1}[\{y\}]))).\]
This is an equivalence $\Pi^1_2$ relation with at least $\omega_2$ classes. As
$E$ remains an equivalence relation after adding a Cohen real we may
apply \cite{HSh} (or \cite{Sh}) to get a perfect set of pairwise
nonequivalent elements. Removing at most one point from this set we find
the desired one.\QED

\noindent{\sc Remark:}\ \ \ In the above theorem, if the ideal $I$
has a $\Pi^1_1$ definition then it is enough to assume that there are
$\omega_1$ respective points $x_\alpha$ for $f$.

\section{An ideal on trees} 
Here we present a Borel $\sigma$--ideal $I^*$ with $\Pi^1_2$ definition which
(in ZFC) satisfies the $\omega_2$--cc but does not have the ccc. Thus if CH
fails then this ideal cannot have the property (M). But we can conclude this
even under CH, if we have enough Cohen reals (see \ref{connotcc}).

\begin{definition}
Let $\Tr$ be the set of all subtrees $T$ of $\fseo$ and let
\[\Trw=\{T\in\Tr: T \mbox{ is well founded }\}.\]
For a tree $T\in\Trw$ let $h_T:T\rightarrow\omega_1$ be the canonical
rank function. For $\alpha<\omega_1$ put
\[A_\alpha=\{T\in\Trw: h_T(\langle\rangle)=\alpha\}.\]
\end{definition}
Clearly $\Tr$ is a closed subset of the product space $2^{({\textstyle
\fseo})}$, so it is a Polish space with the respective topology. The set
$\Trw$ is a $\Pi^1_1$-subset of it and the sets $A_\alpha$ are Borel
subsets of $\Tr$. 

The ideal $I^*$ will live on the space $\Tr$. To define it we introduce
topologies $\tau_\alpha$ on the sets $A_\alpha$ (for $\alpha<\omega_1$). Fix
$\alpha<\omega_1$ for a moment. The basis of the topology $\tau_\alpha$
consists of all sets of the form 
\[U(n,T)=\{T'\in A_\alpha: T'\cap n^{\textstyle \leq n}=T\cap n^{\textstyle
\leq n}\ \&\  h_{T'}\rest n^{\textstyle \leq n}= h_T\rest n^{\textstyle \leq
n}\}\]
for $T\in A_\alpha$, $n\in\omega$. It should be clear that this is a
(countable) basis of a topology. 

\begin{lemma}
\label{baiforalp}
\begin{enumerate}
\item $\langle A_\alpha,\tau_\alpha\rangle$ is a Baire space. 
\item If $\alpha>0$ then there is no isolated point in $\tau_\alpha$.
\item Each Borel (in the standard topology) subset of $A_\alpha$ is
$\tau_\alpha$-Borel. 
\end{enumerate}
\end{lemma}

\Proof 1)\ \ \ Suppose that $O_k\subseteq A_\alpha$ are $\tau_\alpha$-open
dense (for $k<\omega$) and that $T\in A_\alpha$, $n\in\omega$. For each limit
$\beta\leq \alpha$ fix an increasing sequence $\beta(m)\rightarrow\beta$.
Next define inductively sequences $\langle n_k: k\in\omega\rangle$, $\langle
T_k:k\in\omega\rangle$ such that
\begin{description}
\item[$(\alpha)$] $n_k<n_{k+1}<\omega$, $T_k\in A_\alpha$, $n_0=n$, $T_0=T,$
\item[$(\beta)$] $U(n_{k+1},T_{k+1})\subseteq U(n_k,T_k)\cap O_k,$
\item[$(\gamma)$] for each $\nu\in T_k\cap n_k^{\textstyle\leq n_k}$,\\
if $\beta=h_{T_k}(\nu)$ is a limit ordinal then for each $\ell\leq k$ there is
$\nu'\in T_{k+1}\cap n_{k+1}^{\textstyle \leq n_{k+1}}$ such that
$\beta(\ell)\leq h_{T_{k+1}}(\nu')$ and $\nu\vartriangleleft\nu'$,\quad and\\
if $\beta=h_{T_k}(\nu)$ is a successor ordinal then there is $\nu'\in T_{k+1}
\cap n_{k+1}^{\textstyle \leq n_{k+1}}$ extending $\nu$ and such that
$h_{T_{k+1}}(\nu)=\beta-1$.
\end{description}
Put $T^*=\bigcup\limits_{k\in\omega}T_k\cap n_k^{\textstyle \leq n_k}$,
$h^*=\bigcup\limits_{k\in\omega}h_{T_k}\rest n_k^{\textstyle \leq n_k}$. Then
$h^*:T^*\rightarrow\alpha+1$ is a rank function (so $T^*$ is well founded).
Moreover, the condition $(\gamma)$ guarantees that it is the canonical rank
function on $T^*$ and hence $T^*\in\bigcap\limits_{k\in\omega}O_k\cap U(n,T)$.

\noindent The assertions 2) and 3) should be clear. \QED

\begin{definition}
\label{theideal}
The ideal $I^*$ consists of subsets of Borel sets $B\subseteq\Tr$ such that
\[(\forall 0<\alpha<\omega_1)(B\cap A_\alpha\ \mbox{ is }
\tau_\alpha\mbox{-meager}).\]
\end{definition}

\begin{proposition}
\label{notcc}
$I^*$ is a non-trivial Borel $\sigma$--ideal on $\Tr$ which does not satisfy
the ccc but satisfies the $\omega_2$--cc.
\end{proposition}

\Proof As $\{A_\alpha:\alpha<\omega_1\}$ is a partition of $\Trw$ and each
$A_\alpha$ is not $\tau_\alpha$-meager (by \ref{baiforalp}), we get that $I^*$
is a non-trivial $\sigma$-ideal on $\Tr$. Since each set $A_\alpha$ is not in
$I^*$, the ccc fails for $I^*$. If $\{B_\zeta:\zeta\in\omega_2\}$ is a family
of $I^*$-almost disjoint $I^*$--positive Borel sets then for some $\alpha<
\omega_1$ and a set $Z\in [\omega_2]^{\textstyle \omega_2}$, $\{B_\zeta\cap
A_\alpha:\zeta\in Z\}$ is a family of $\tau_\alpha$-Borel
$\tau_\alpha$--non-meager ($\tau_\alpha$-meager)-almost disjoint sets. This
contradicts the fact that the topology $\tau_\alpha$ has a countable basis
(and each basic set is not $\tau_\alpha$-meager). \QED
\medskip

Now, we want to estimate the complexity of $I^*$.
\begin{proposition}
\label{complexity}
The ideal $I^*$ has a $\Pi^1_2$ definition.
\end{proposition}

\Proof Fix $\alpha<\omega_1$. A basic open subset of $A_\alpha$ (in the
topology $\tau_\alpha$) is determined by a pair $\langle F,h\rangle$, where
$F\subseteq\fseo$ is a finite tree and $h:F\rightarrow\alpha+1$ is a rank
function such that $h(\langle\rangle)=\alpha$. The basic sets $U(F,h)$ are
Borel subsets of $A_\alpha$ (in the standard topology). Let
\[\begin{array}{ll}
\Phi_{\nwd}(a,\alpha)\equiv &\mbox{``}a\mbox{ is a Borel code }\ \&\\
\ &(\forall \langle F,h\rangle )(\exists \langle F',h'\rangle )[F\subseteq
F'\ \& \ h\subseteq h'\ \&\ U(F',h')\cap\# a=\emptyset]\mbox{''.}\\
\end{array}\]
It should be clear that $\Phi_{\nwd}(\cdot,\alpha)$ is a $\Pi^1_1$-formula.
Consequently, the formula ``$a$ is a Borel code and the set $\# a\cap
A_\alpha$ is $\tau_\alpha$-nowhere dense'' is $\Pi^1_1$. Let
\[\begin{array}{ll}
\Phi_{\meager}(a,\alpha)\equiv &\mbox{``} a\mbox{ is a Borel code }\ \&\\
\ &(\exists\langle c_n:n\in\omega\rangle)[(\forall n)\Phi_{\nwd}(c_n,\alpha)\ 
\&\ \# a\cap A_\alpha\subseteq\bigcup\limits_{n\in\omega}\# c_n]\mbox{''.}\\
\end{array}\]
This formula says that the intersection $\# a\cap A_\alpha$ is
$\tau_\alpha$-meager; it is $\Sigma^1_2$. On the other hand, its negation is
equivalent to
\[\begin{array}{l}
\Phi_{\nmeager}(a,\alpha)\equiv \mbox{``}a\mbox{ is not a Borel code }
\mbox{ or}\\
\qquad (\exists\langle F,h\rangle)(\exists\langle c_n:n\in\omega\rangle)[
(\forall n)\Phi_{\nwd}(c_n,\alpha)\ \&\ U(F,h)\subseteq\# a\cup
\bigcup\limits_{n\in\omega}\#c_n]\mbox{''}
\end{array}\] 
which is $\Sigma^1_2$ too. Consequently, 
\smallskip
\begin{center}
``$a$ is a Borel code such that $\#a\cap A_\alpha$ is $\tau_\alpha$-meager'' 
\end{center}
\smallskip

\noindent is absolute between all transitive models of (a large enough part
of) ZFC containing $a,\alpha,\ldots$. Now, we can give a $\Pi^1_2$-definition
of the ideal $I^*$: 
\[\begin{array}{ll}
\Phi(a)\equiv& a \mbox{ is a Borel code }\ \&\\
\ &(\forall E\subseteq\omega\times\omega)[E \mbox{ is not well founded or }
\langle \omega,E\rangle \not\models\Theta^*\mbox{ or }\\
\ &a\mbox{ is not encoded in }\langle \omega,E\rangle \mbox{ or }\langle
\omega,E\rangle \not\models\V=\L[a] \mbox{ or }\\
\ &\langle \omega,E\rangle \models \mbox{``}(\forall\alpha\in\omega_1)
(\Phi_{\meager}(a,\alpha))]\mbox{''.}\\
\end{array}\]
The $\Theta^*$ above is a sentence carrying a large part of information on
ZFC, in particular for each transitive set $M$, $M\models\Theta^*$ should
imply that $M$ is $A$-adequate for all $A\in M$ and the absoluteness of
$\Sigma^1_1$, $\Pi^1_1$ formulas holds for $M$ (see \cite[Ch. 2, \S15]{Je}).
Thus, if $E$ is a well founded relation on $\omega$, $a$ is a real encoded in
$E$ (which means that $a$ belongs to the transitive collapse of $E$) and
$\langle \omega,E\rangle \models \Theta^* + \V=\L[a]$ then the transitive
collapse of $E$ is $\L_\alpha[a]$ for some $\alpha<\omega_1$. Now, if
$\Phi(a)$ holds then for all $\alpha<\omega_1$ such that $\L_\alpha[a]$ is
$a$-adequate we have 
\[\L_\alpha[a]\models \mbox{``}(\forall\beta<\omega_1)\Phi_{\meager}(a,\beta)
\mbox{''}.\]
The absoluteness implies that if $\L_\alpha[a]\models\Phi_{\meager}(a,\beta)$
then $\Phi_{\meager}(a,\beta)$. Consequently, $\Phi(a)$ implies
$(\forall\beta<\omega_1)\Phi_{\meager}(a,\beta)$. Suppose now that $a$ is a
Borel code for which $\Phi(a)$ fails. Then for some $\alpha<\omega_1$
(suitably closed) we have
\[\L_\alpha[a]\models(\exists\beta)(\neg\Phi_{\meager}(a,\beta))\]
and hence for some $\beta<\alpha$ we get
\[\L_\alpha[a]\models\Phi_{\nmeager}(a,\beta).\]
Once again the absoluteness implies that $\Phi_{\nmeager}(a,\beta)$ holds.
Consequently, $\neg\Phi(a)$ implies that there is $\beta<\omega_1$ such that
$\Phi_{\nmeager}(a,\beta)$ and finally, 
\[\Phi(a)\iff a\mbox{ is a Borel code }\ \&\ \# a\in I^*.\]
It should be clear that $\Phi(x)$ is a $\Pi^1_2$ formula. \QED

\begin{conclusion}
\label{connotcc}
Assume that either CH fails or every $\Delta^1_2$ set of reals has the Baire
property. Then the ideal $I^*$ does not have the property (M) [and: it
satisfies the $\omega_2$-cc, it does not satisfy the ccc and has a $\Pi^1_2$
definition].  
\end{conclusion}

\Proof If we are in the situation of $\neg$CH then (M) cannot hold because of
$\omega_2$-cc. So assume that all $\Delta^1_2$ sets of reals have Baire
property.\\ 
Suppose that $f:\can\rightarrow\can$ is a Borel function (of course, it is
coded by a real). Let $\langle c_\alpha:\alpha<\omega_2\rangle$ be a
$\C_{\omega_2}$-generic sequence of Cohen reals. In
$\V[c_\alpha:\alpha<\omega_2]$ we have that for sufficiently large $\alpha$,
$f^{-1}[\{c_\alpha\}]\in I^*$. Since the definition of $I^*$ is absolute (and
it involves no parameters), 
\begin{center}
$\V[c_\alpha:\alpha<\omega_2]\models f^{-1}[\{c_\beta\}]\in I^*$\quad implies
\quad $\L[f][c_\beta]\models f^{-1}[\{c_\beta\}]\in I^*$. 
\end{center}
By density arguments and the ccc of $\C_{\omega_2}$ we conclude that
\[\L[f]\models\mbox{``}\forces_{\C}f^{-1}[\{\dot{c}\}]\in I^*\mbox{''}\]
where $\dot{c}$ is the canonical $\C$-name for the generic Cohen real. Our
assumption that $\Delta^1_2$ sets of reals have the Baire property is
equivalent to the statement 
\begin{center}
``for each real $r$ there exists a Cohen real over $\L[r]$'' 
\end{center}
(see e.g. \cite{JuSh}). Let $c\in\V$ be a Cohen real over
$\L[f]$. Then 
\[\L[f][c]\models f^{-1}[\{c\}]\in I^*.\]
By Shoenfield absoluteness we conclude $\V\models f^{-1}[\{c\}]\in I^*$.
Consequently, $f$ cannot witness the property (M) for $I^*$.\QED

\section{Other variants of (M)}
The following natural question concerning the hereditary behavior of our
conditions (B), (M) and (D) arises: 
\begin{quotation}
\noindent If an ideal $I$ on $X$ satisfies one of those conditions and
$E\notin I$ is a Borel subset of $X$, is it true that $I\cap P(E)$ satisfies
the same condition (on $E$; at this moment we consider the ideal on the space
$E$ and modify respectively the sense of (B), (M) and (D)))?
\end{quotation}
When we do not assume a group structure on $X$ and the invariance of $I$, one
can easily construct $I$ which makes the answer negative. If the invariance of
$I$ is supposed, the problem becomes less trivial.
\medskip

\noindent{\sc Example:}\ \ \ Let $I$ be the ideal of null sets
with respect to $1$-dimensional Hausdorff measure (cf. \cite{F}) on
${\bf R}^{2}$ treated as an additive group. Clearly, $I$ is invariant. It
satisfies (M) since the fibers of the continuous function given by
$f(x,y)=x$ are the lines $\{ x\} \times {\bf R}$, for $x\in {\bf
R}$, which are not in $I$. However, $1$-dimensional Hausdorff
measure on ${\bf R}$ coincides with the linear Lebesgue measure,
so $I\cap P(\{ x\}\times {\bf R})$ does not satisfy (B) since ccc works there.
\medskip

Further, we shall concentrate on some hereditary versions of (M)
connected with open sets.

\begin{definition}
Assume that $I$ is an ideal on an uncountable Polish space $X$ and that each
open nonempty subset of $X$ is not in $I$.
\begin{enumerate}
\item We say that $I$ has {\em property (M) relatively to a set $E\subseteq
X$} (in short, $I$ has (M) rel $E$ ) if and only if\\
there is a Borel function $f:E\rightarrow X$ whose all fibers $f^{-1}[\{ x\}
]$, $x\in X$, are not in $I$.
\item We say that $I$ has {\em property (M')} if and only if\\
it has property (M)  rel $U$ for each nonempty open set $U\subseteq X$.
\item Let $x\in X$. We say that $I$ has {\em property $(M'_{x})$} if and only
if\\
it has property (M) rel $U$ for each open neighborhood $U$ of $x$.
\item We say that $I$ has {\em property $(M^{*})$} if and only if\\
there is a Borel function $f:X\rightarrow X$ such that $f^{-1}[\{x\}]\cap
U\notin I$ for any $x\in X$ and open non-void $U\subseteq X$.
\end{enumerate}
\end{definition}
\smallskip

\noindent{\sc Remarks:}
\begin{enumerate}
\item Plainly, (M) in the sense of Definition 1.1 is (M) rel $X$. Moreover,
if $E$ is Borel, (M) rel $E$ implies (M) rel $X$ since one can consider any
Borel extension of the respective Borel function $f:E\rightarrow X$ to the
whole $X$.
\item Obviously, $I$ has $(M')$ iff it has $(M'_{x})$ for each $x$ from a fixed
dense set in $X$. In some cases $(M'_{x})$ satisfied for one point $x$ implies
$(M')$. This holds if $x$ is taken from a fixed dense set whose any two points
$s,t$ have homeomorphic open neighborhoods $U_{s},U_{t}$ and the homeomorphism
preserves $I$. In particular, if $X$ is a group, $x\in X$ is fixed, $Q$ is a
dense set in $X$ and $I$ is $Q$-invariant (i.e. $E\in I$ and $t\in Q$ imply
$E+t\in I$) then $(M'_{x})$ guarantees $(M')$.
\item Studies of $(M^{*})$ were initiated in \cite{Ba1} where some examples are
given. Condition $(M^{*})$ for the ideal of nowhere dense sets means the
existence of a Borel function $f$ from $X$ onto $X$ with dense fibers. For
$X=(0,1)$ such functions are known as being strongly Darboux (cf.
\cite{Br}). From Mauldin's proof in \cite{Md*} it follows that the
$\sigma$-ideal generated by closed Lebesgue null sets satisfies $(M^{*})$
(for details, see \cite{Ba1}).
\item Evidently, $(M^{*})\Rightarrow (M')\Rightarrow (M)$. It is interesting
to know whether those implications can be reversed. A simple method
producing ideals ($\sigma$-ideals) with property (M) and without (M')
follows from the example given in \cite{Ba1}. If one part of an ideal $I$,
defined on a Borel set $B\subseteq X$ has (M) rel $B$, and the remaining
part, defined on $X\setminus B$ with $int(X\setminus B)\neq \emptyset$ has
not (M) rel $X\setminus B$ (for instance, ccc holds there) then the
ideal is as desired. If we want $I$ to be a $\sigma$-ideal invariant
in the group $X$ then the situation is different.
\end{enumerate}

\begin{theorem}
\label{versions}
Suppose that $I$ is an invariant $\sigma$-ideal on {\bf R}. If the ideal $I$
has the property (M) then it has the property $(M^*)$. 
\end{theorem}

\Proof We start with two claims of a general character.
\smallskip
\begin{claim}
\label{cl2}
Suppose that $f:{\bf R}\longrightarrow{\bf R}$ is a Borel function. Then there
exists a perfect set $P\subseteq{\bf R}$ such that $f^{-1}[P]$ is both meager
and Lebesgue null.
\end{claim}

\noindent {\em Proof of the claim:}\ \ \ It follows directly from the fact
that each perfect set can be divided into continuum disjoint perfect sets and
both the ideal of meager sets and the ideal of null sets satisfy ccc.

\begin{claim}
\label{cl3}
Suppose that $U\subseteq{\bf R}$ is an open set and $A\subseteq{\bf R}$ is a
meager null set. Then there exist disjoint sets $D_k\subseteq U$ and reals
$x_k$ (for $k\in\omega$) such that $\bigcup\limits_{k\in\omega} D_k$ is
nowhere dense and $A\subseteq\bigcup\limits_{k\in\omega} D_k+x_k$.
\end{claim}

\noindent{\em Proof of the claim:}\ \ \ First we prove the claim under the
assumption that $A$ is nowhere dense (and null). For this choose an interval
$J\subseteq U$. Since $A$ is null we can find (open) intervals $\langle J_k:
k\in\omega\rangle$ such that
\[A\subseteq\bigcup_{k\in\omega} J_k\quad\quad\mbox{ and }\quad\quad
\sum_{k\in\omega} |J_k|<|J|\]
(here $|J|$ stands for the length of $J$).
Choose disjoint open intervals $J_k^*\subseteq J$ (for $k\in\omega$) such
that $|J^*_k|=|J_k|$ (for $k\in\omega$). Let reals $x_k$ be such that
$J_k^*+x_k=J_k$. For $k\in\omega$ we put
\[D_k=(A\cap J_k)-x_k.\]
Clearly $A\subseteq\bigcup\limits_{k\in\omega} D_k+x_k$. Since the sets
$D_k$ are nowhere dense and contained in disjoint open intervals we get that
their union $\bigcup\limits_{k\in\omega}D_k$ is nowhere dense.

If now $A$ is just a meager null set then we represent it as a union
$A=\bigcup\limits_{n\in\omega}A_n$, where each set $A_n$ is nowhere dense
and null. Take disjoint open sets $U_n\subseteq U$ (for $n\in\omega$) and
apply the previous procedure to each pair $\langle A_n,U_n\rangle$ getting
suitable $D^n_k,x^n_k$. Then $\langle D^n_k,x^n_k: n,k\in\omega\rangle$ is as
required for $A$ and $U$ (proving the claim).
\medskip

Suppose now that $I$ is an invariant $\sigma$-ideal on ${\bf R}$ and that
$f:{\bf R}\longrightarrow{\bf R}$ is a Borel function witnessing the
property $(M)$ for $I$. By Claim \ref{cl2} we find a perfect set
$P\subseteq{\bf R}$ such that $f^{-1}[P]$ is both meager and null.

\noindent Let $\{U_n: n\in\omega\}$ be an enumeration of all rational open
intervals in {\bf R}. As finite union of nowhere dense sets is nowhere dense,
we may apply Claim \ref{cl3} and choose inductively (by induction on
$n\in\omega$) closed sets $D_{n,k}$ and reals $x_{n,k}$ such that
\begin{description}
\item[(a)\ \ ] $D_{n,k}$ (for $n,k\in\omega$) are pairwise disjoint
\end{description}
and for each $n\in\omega$:
\begin{description}
\item[(b)\ \ ] $\bigcup\limits_{k\in\omega}D_{n,k}$ is a nowhere dense subset
of $U_n$,
\item[(c)\ \ ] $f^{-1}[P]\subseteq\bigcup\limits_{k\in\omega}D_{n,k}+x_{n,k}$.
\end{description}
Now we define a function $f^*:{\bf R}\longrightarrow P$ by
\[f^*(x)=\left\{
\begin{array}{ll}
f(x+x_{n,k}) &\mbox{ if } x\in D_{n,k},\quad n,k\in\omega,\quad \mbox{and }
f(x+x_{n,k})\in P,\\
y_0          &\mbox{ otherwise,}
\end{array}
\right.\]
where $y_0$ is a fixed element of $P$. Clearly the function $f^*$ is Borel.
We claim that it witnesses the property $(M^*)$ for $I$. Why? Suppose that
$U\subseteq{\bf R}$ is an open nonempty set and $y\in P$. Take $n\in\omega$
such that $U_n\subseteq U$. We know that $f^{-1}[P]\subseteq
\bigcup\limits_{k\in \omega}D_{n,k}+x_{n,k}$ and $f^{-1}[\{y\}]\notin I$.
Since $I$ is $\sigma$-complete, we find $k\in\omega$ such that
\[f^{-1}[\{y\}]\cap (D_{n,k}+x_{n,k})\notin I.\]
As $I$ is translation invariant we get
\[(f^{-1}[\{y\}]-x_{n,k})\cap D_{n,k}\notin I.\]
Clearly $(f^{-1}[\{y\}]-x_{n,k})\cap D_{n,k}\subseteq (f^*)^{-1}[\{y\}]\cap
U$ so the last set does not belong to $I$. The theorem is proved. \QED
\medskip

\noindent In the above proof the use of Lebesgue measure is important.
In Theorem 5.2 we can replace {\bf R} by e.g. $\can$ but we
do not know if we can have a corresponding result for all Polish groups.
Moreover we do not know if we can omit the assumption of
$\sigma$-completeness of $I$ (i.e. prove the theorem for (finitely additive)
ideals on {\bf R} which do not contain nonempty open sets).

\noindent For the question of dependences between $(M)$, $(M')$
and $(M^*)$ the following simple theorem seems useful.

\begin{theorem}
\label{usefull}
Let $I$ be an ideal on $X$ and $\{ U_{n}: n\in \omega \}$ - a base of open
sets in $X$. For $f:X\rightarrow X$ define $H_{n}(f)=\{ y\in X: U_{n}\cap
f^{-1}[\{ y\} ]\notin I\}$.
\begin{enumerate}
\item Let $x\in X$. Condition $(M'_{x})$ holds iff there is a Borel function
$f:X\rightarrow X$ such that, for each $n\in \omega$ with $x\in U_{n}$, the
set $H_{n}(f)$ contains a perfect set.
\item Condition $(M')$ holds iff there is a Borel function $f:X\rightarrow X$
such that each set $H_{n}(f)$, $n\in \omega$, contains a perfect set.
\item Condition $(M^{*})$ holds iff there is a Borel function $f:X\rightarrow
X$ such that $\bigcap\limits_{n\in \omega} H_{n}(f)$ contains a perfect set. 
\item Suppose that $I$ is a $\sigma$-ideal, $f:X\rightarrow X$ is a Borel
function such that $f^{-1}[\{x\}]\notin I$ for all $x\in X$, and each set
$H_{n}(f)$, $n\in \omega$, either is countable or contains a perfect set.
Then $(M'_{x})$ holds for some $x\in X$.
\end{enumerate}
\end{theorem}

\Proof $1.$ Necessity. Fix $n\in \omega $ such that $x\in U_{n}$. Then
(M) rel $U_{n}$ holds. It suffices to extend the respective Borel function
$g:U_{n}\rightarrow X$ to a Borel function defined on the whole $X$.\\
Sufficiency. Fix $n\in \omega$ such that $x\in U_{n}$. Let $P\subseteq
H_{n}(f)$ be a perfect set and let $B=f^{-1}[P]$. Extend $f\rest B$ to a
Borel function $g:U_{n}\rightarrow P$. Consider a Borel function $h$ from
$P$ onto $X$. Then $h\circ g$ witnesses (M) rel $U_{n}$.\\
The proofs of $2$ and $3$ are analogous.

$4.$ Suppose it is not the case. Thus for each $x\in X$ choose $n_{x}\in
\omega$ such that $I$ has not (M) rel $U_{n_{x}}$. For $T=\{n_{x}:x\in X\}$ we
get $X=\bigcup\limits_{n\in T}U_n$. Since $I$ is a $\sigma$--ideal, therefore,
by the properties of $f$, we get $X=\bigcup\limits_{n\in T}H_{n}(f)$. Hence,
by the assumption, there is $n\in T$ such that $H_{n}(f)$ contains a perfect
set.  Now, as in the proof of $1$, we infer that $I$ has (M) rel $U_{n}$, a
contradiction.\QED 

\noindent{\sc Example:}\ \ \ Let $I_2$ ($I_\omega$) be a Mycielski ideal
(defined in \cite{My3}) on $\can$ ($\baire$, respectively) generated by the
respective system $\{K_t:t\in 2^{<\omega}\}$ of infinite sets in $\omega$. It
is shown in \cite{BaRo} that $I_2$, $I_\omega$ satisfy (M). By
theorem~\ref{versions} the ideal $I_2$ satisfies $(M^*)$. However, it is not 
clear if $I_\omega$ does. Let us modify the proof of $(M)$ for $I_\omega$
to get (M) rel $[s]$ where $[s]$, for $s\in \omega^{<\omega }$, is a basic
open set in $\baire$. To this end choose $K_t$ such that $K_t\cap lh(s)=
\emptyset$. Assume that $\omega\setminus (K_t\cup lh(s))=\{n_0,n_1,\ldots\}$
and define $f:[s]\rightarrow [s]$ by 
\[f(s\conc\langle x_0,x_1,\ldots\rangle)=s\conc\langle x_{n_0},x_{n_1},\ldots
\rangle.\]
As in \cite{BaRo} we observe that $f$ realizes (M) rel $[s]$. Consequently,
$I_\omega$ satisfies (M').

\begin{problem}
\begin{enumerate}
\item Does Mycielski ideal $I_\omega$ on $\baire$ satisfy $(M^{*})$ ?
\item Does there exist a translation invariant $\sigma$-ideal $I$ on
$\baire\equiv {\bf Z}^\omega$ which satisfies $(M)$ but not $(M^*)$ ?
\item Is there (necessarily finitely additive) invariant ideal $I$ on {\bf R}
with $(M)$ but without $(M^*)$ (and such that no nonempty open set is in $I$)?
\end{enumerate}
\end{problem}

\noindent{\sc Remark:}\ \ \ By assertion $4$ of Theorem~\ref{usefull}, if
the axiom of determinacy AD (cf. \cite{Je}) is assumed, we get $(M)\iff
(\exists x\in X) (M'_{x})$ for any $\sigma$-ideal and consequently, $(M)
\iff (M')$ for any invariant $\sigma $-ideal.
The following operation (cf. \cite{BaRo}) $\Phi _{I}:P(X\times X)\rightarrow
P(X)$ plays an important role. Namely, for $E\subseteq X\times X$, let
$\Phi _{I}(E)=\{ y\in X: E^{y}\notin I\}$ where $E^{y}=\{ x\in X: \langle x,y
\rangle \in E\}$. If $\Phi_{I}$ sends Borel sets into analytic sets
(respectively, into projective sets when projective determinacy is assumed)
then statement $4$ of Theorem~\ref{usefull} works.

\end{document}